\documentclass{article}

\usepackage{amssymb, latexsym}

\usepackage{graphicx}

\usepackage{amsmath}

\newtheorem{theorem}{Theorem}[section]

\newtheorem{corollary}[theorem]{Corollary}

\newtheorem{definition}[theorem]{Definition}

\newtheorem{lemma}[theorem]{Lemma}

\newtheorem{proposition}[theorem]{Proposition}

\newtheorem{remark}[theorem]{Remark}

\begin{document}
\title{\Large\bf  Hamilton's Harnack inequality and the $W$-entropy formula on complete Riemannian manifolds}
\author{\ \  Xiang-Dong Li
\thanks{Research supported by NSFC No. 11371351, Key Laboratory RCSDS, CAS, No. 2008DP173182, and a
Hundred Talents Project of AMSS,
CAS.}
}

\maketitle

\begin{minipage}{120mm}
{\bf Abstract}. In this paper, we prove Hamilton's Harnack inequality and the gradient estimates of the logarithmic heat kernel for the Witten Laplacian on complete Riemainnian manifolds. As applications, we prove the $W$-entropy formula for the Witten Laplacian on complete Riemannian manifolds, and prove a family of logarithmic Sobolev inequalities on complete Riemannian manifolds with natural geometric condition.
\end{minipage}

\section{Introduction}

Differential Harnack inequality is an important tool in the study of heat equation on complete Riemannian manifolds \cite{LY, Da, Li05, SY, CLN, Chow}. In \cite{H1}, Hamilton proved a Harnack inequality for  positive and bounded solutions to the heat equation on compact Riemannian manifolds. More precisely, let $M$ be a complete Riemannian manifold, and suppose that  there exists a constant $K\geq 0$ such that
\begin{eqnarray*}
Ric\geq -K.
\end{eqnarray*}
Then, for any positive and bounded solution $u$ to the heat equation
\begin{eqnarray*}
\partial_t u=\Delta u,
\end{eqnarray*}
it holds
\begin{eqnarray*}
|\nabla \log u|^2\leq \left({1\over t}+2K\right)\log(A/u),\ \ \
\forall x\in M, t>0,
\end{eqnarray*}
where
\begin{eqnarray*}
A:=\sup\limits\{u(t, x): x\in M, t\geq 0\}.
\end{eqnarray*}

In \cite{Ko}, Kotschwar extended Harmilton's Harnack inequality to
complete Riemannian manifolds with non-negative Ricci curvature. As
application, he proved the following gradient estimates of the
logarithmic of the heat kernel $p_t(x, y)$ of the Laplace-Beltrami
operator on complete Riemannian manifolds with non-negative Ricci
curvature. More precisely, for all $\delta>0$, there exists a
constant $C=C(n, \delta)>0$, such that
\begin{eqnarray*}
|\nabla \log p_t(x, y)|^2\leq {2\over t}\left(C+{d^2(x, y)\over (4-\delta)t}\right), \ \ \ \forall x, y\in M, t>0.
\end{eqnarray*}

The purpose of this paper is to prove Hamilton's Harnack inequality and gradient estimates of the logarithmic of the heat kernel of the Witten Laplacian  on complete Riemannian manifolds with weighted volume measure.  As application, we derive the $W$-entropy formula  and prove the rigidity theorem of the $W$-entropy for the Witten Laplacian on complete Riemannian manifolds with natural Bakry-Emery curvature condition. Indeed, the relation between the $W$-entropy functional and the differential Harnack inequality has been recently investigated by many authors \cite{P1, N1, N2, Chow, Ec, LNVV}.

Let $M$ be a complete Riemannian manifold with a fixed Riemnnian metric $g$, $\phi\in C^2(M)$ and $d\mu=e^{-\phi}dv$, where $v$ is the Riemannian volume measure on $(M, g)$, and let
\begin{eqnarray*}
L =\Delta -\nabla \phi\cdot\nabla \label{WL}
\end{eqnarray*}
be the Witten Laplacian on $(M, g)$ with respect to the weighted volume measure $\mu$. For all $u, v\in C^\infty_0(M)$, the following integration by parts formula holds
\begin{eqnarray*}
\int_M \langle \nabla u, \nabla v\rangle d\mu=-\int_M Lu vd\mu=-\int_M uLvd\mu.
\end{eqnarray*} In \cite{BE}, Bakry and Emery proved that for all $u\in C_0^\infty(M)$,
\begin{eqnarray}
L|\nabla u|^2-2\langle \nabla u, \nabla L u\rangle=2|\nabla^2
u|^2+2Ric(L)(\nabla u, \nabla u),\label{BWF}
\end{eqnarray}
where
$$Ric(L)=Ric+\nabla^2\phi.$$ The formula $(\ref{BWF})$
can be viewed as a natural extension of the Bochner-Weitzenb\"ock
formula. The quantity $Ric(L)=Ric+\nabla^2 \phi$, called the
Bakry-Emery Ricci curvature on the weighted Riemannian manifolds
$(M, g, \phi)$, plays as a good substitute of the Ricci curvature in
many problems in comparison geometry and analysis on complete
Riemannian manifolds with smooth weighted volume measures. See
\cite{BE, FLZ, FLL, Li05, Li11a, Li11b, Q, WW, MW, MW2} and
reference therein.

Following \cite{BE, Lot, Li05}, we introduce the $m$-dimensional Bakry-Emery Ricci curvature on $(M, g, \phi)$ by
\begin{eqnarray*}
Ric_{m, n}(L):=Ric+\nabla^2\phi-{\nabla\phi\otimes \nabla \phi\over m-n},
\end{eqnarray*}
where $m\geq n$ is a constant, and $m=n$ if and only if $\phi$ is a constant. When $m=\infty$, we have $Ric_{\infty, n}(L)=Ric(L)$.
When $m\in \mathbb{N}$, the $m$-dimensional Bakry-Emery Ricci curvature $Ric_{m, n}(L)$ has a very natural geometric interpretation.
Indeed, consider the warped product metric on $M^n\times S^{m-n}$ defined by
\begin{eqnarray*}
\widetilde{g}=g_M\bigoplus e^{-{2\phi\over m-n}}g_{S^{m-n}}.\label{WPM}
\end{eqnarray*}
where $S^{m-n}$ is the unit sphere in $\mathbb{R}^{m-n+1}$ with the standard metric $g_{S^{m-n}}$. By a classical result in Riemaniann geometry, the quantity $Ric_{m, n}(L)$ is equal to the Ricci curvature of the above warped product metric $\widetilde g$ on $M^n\times S^{m-n}$ along the horizontal vector fields. See \cite{Lot, Li05, WW}.

We now state the main results of this paper. The first result of this paper is the following improved version of the Harnack Harnack inequality for the heat equation of the Witten Laplacian on complete Riemannian manifolds.

\begin{theorem}\label{Thm0} Let $M$ be a complete Riemannian manifold, and $\phi\in C^2(M)$. Suppose that there exists a constant $K\geq 0$ such that
\begin{eqnarray*}
Ric(L)\geq -K.
\end{eqnarray*}
Let $u$ be a positive and bounded solution to the heat equation
\begin{eqnarray*}
\partial_t u=L u.
\end{eqnarray*}
Let $A:=\sup\limits\{u(t, x): x\in M, t\geq 0\}$. Then
\begin{eqnarray}
|\nabla \log u|^2&\leq& {2K\over 1-e^{-2Kt}}\log (A/u). \label{BLH}
\end{eqnarray}
\end{theorem}

The following result extends the Hamilton Harnack inequality to the Witten Laplacian on complete Riemannian manifolds.

\begin{theorem}\label{Thm02}  Under the same condition as in Theorem \ref{Thm0}, the Halmilton  Harnack inequality holds
\begin{eqnarray}
|\nabla \log u|^2\leq \left({1\over t}+2K\right)\log(A/u).\label{Ham}
\end{eqnarray}
\end{theorem}

Letting $t\rightarrow \infty$ in Theorem \ref{Thm02}, we can derive the following gradient estimate and the Liouville theorem due to Brighton \cite{Br}. See also Remark \ref{rem4} below.

 \begin{corollary}  \label{cor1} Let $u$ be a bounded and positive $L$-harmonic function $u$ on $M$. Then, under the same condition as in Theorem \ref{Thm0}, we have
\begin{eqnarray*}
|\nabla \log u|\leq 2K\log(A/(u-\inf\limits u)),
\end{eqnarray*}
where $A=\sup\limits u-\inf\limits u$. In particular, if $Ric+\nabla^2\phi\geq 0$, then every bounded $L$-harmonic function must be constant.
\end{corollary}

Let $p_t(x, y)$ be the heat kernel of $L$ with respect to the reference measure $\mu$, i.e., the fundamental solution of the heat equation
$\partial_t u=Lu$. The following result extends Hamilton and Kotschwar's gradient estimate to the logarithmic heat kernel of the Witten Laplacian on complete Riemannian manifolds with natural geometric condition.

\begin{theorem}\label{Thm1} Let $M$ be a complete Riemannian manifolds, and $\phi\in C^2(M)$. Suppose that there exist some constant $m\geq n$ and $K\geq 0$ such that
$$Ric_{m, n}(L)\geq -K.$$
Then, for all $T>0$, there exists a constant $C(K, m, n, T)>0$ such that
\begin{eqnarray*}
|\nabla \log p_t(x, y)|\leq C(K, m, n, T)\left({d(x, y)\over t}+{1\over \sqrt{t}}\right)
\end{eqnarray*}
holds for all $x, y\in M$ and $t\in (0, T]$.
\end{theorem}

The following result gives the higher order gradient estimates for the logarithmic heat kernel of the Witten Laplacian on complete Riemannian manifolds with natural geometric condition.

\begin{theorem}\label{Thm2}  Let $(M, g)$ be an $n$-dimensional complete Riemannian manifold. Suppose that $\|\nabla^k{\rm Riem}\|\leq C_k$ for some constant $C_k>0$, $0\leq k\leq N$,
$\phi\in C^{N+1}(M)$ with $\nabla\phi \in C_b^{N}(M)$, and there exists some $m\geq n$, $K\geq 0$ such that \ $Ric_{m, n}(L)\geq -K$.
Then, for all $T>0$, there exists a constant $C_N(K, m, n, T)>0$ such that
\begin{eqnarray*}
|\nabla^N \log p_t(x, y)|\leq C_N(K, m, n, T)\left({d(x, y)\over t}+{1\over \sqrt{t}}\right)^N
\end{eqnarray*}
holds for all $x, y\in M$ and $t\in (0, T]$.
\end{theorem}

As an application of Theorem \ref{Thm1} and Theorem \ref{Thm2}, we
prove the following $W$-entropy formula for the Witten Laplacian on
complete Riemannian manifolds, which extends a previous result due
to Ni \cite{N1, N2} for the case $m=n$ and $L=\Delta$.

\begin{theorem}\label{Thm3}  Let $(M, g)$ be a complete Riemannian manifold. Suppose that $\|\nabla^k{\rm Riem}\|\leq C_k$ for some constant $C_k>0$, $0\leq k\leq 2$,
$\phi\in C^{4}(M)$ with $\nabla\phi \in C_b^{3}(M)$, and there
exists some $m\geq n$, $K\geq 0$ such that \ $Ric_{m, n}(L)\geq -K$.
Let $u(t, x)={e^{-f}\over (4\pi t)^{m/2}}$
 be the heat kernel of the heat equation
${\partial u\over \partial t}=Lu$. Let
\begin{eqnarray*}
W_m(u, t)=\int_M \left(t|\nabla f|^2+f-m\right)
{e^{-f}\over (4\pi t)^{m/2}}d\mu \label{WW1}
\end{eqnarray*}
be the $W$-entropy of the Witten Laplacian $L$ on $(M, g, \mu)$.
Then, for all $t>0$, we have
\begin{eqnarray}
{d W_m(u, t)\over dt}&=&-2\int_M
t \left(\left|\nabla^2 f-{g\over 2\tau}\right|^2+Ric_{m,
n}(L)(\nabla f,
\nabla f)\right)ud\mu\nonumber\\
& & \ \ -{2\over m-n}\int_M t \left({\nabla\phi\cdot\nabla
f}+{m-n\over 2t}\right)^2ud\mu.\label{WW2}
\end{eqnarray}
In particular, if $Ric_{m, n}(L)\geq 0$, then $W_m(u, t)$ is
non-increasing in $t$, and $W_m(u, t)$ attains its minimum at some
point $t_0>0$, i.e., ${dW_m(u, t)\over dt}=0$ holds at $t=t_0$, if
and only if $(M, g)$ is isometric to Euclidean space $\mathbb{R}^n$,
$m=n$, $\phi\equiv C$ for a constant $C\in \mathbb{R}$, and
\begin{eqnarray*}
u(t, x)={e^{-{\|x\|^2\over 4t}}\over (4\pi t)^{n/2}}, \ \ \ \ \forall x\in \mathbb{R}^n, t>0.
\end{eqnarray*}
\end{theorem}

In our previous paper \cite{Li11a}, we proved the $W$-entropy
formula $(\ref{WW2})$ for the Witten Laplacian on complete
Riemannian manifolds satisfying the condition in Theorem \ref{Thm3}
and the following uniform volume lower bound condition: there exists
a constant $r_0>0$ such that
\begin{eqnarray}
\mu_0=\inf\limits\{\mu(B(x, r_0): x\in M\}>0.\label{mulow}
\end{eqnarray}
By Theorem $2.8$ in  \cite{Li10}, on complete Riemannian manifolds
with $Ric_{m, n}(L)\geq -K$ and satisfying the uniformly lower bound
condition $(\ref{mulow})$, the following Sobolev inequality holds
\begin{eqnarray*}
\|f\|^2_{2m\over m-2}\leq C\left(\int_M |\nabla f|^2d\mu+\int_M
|f|^2d\mu\right), \ \ \forall f\in C_0^\infty(M).
\end{eqnarray*}
Thus, at least from a geometric and analytic point of view, the
uniformly volume lower bound condition $(\ref{mulow})$ is too
restrictive for the validity of the $W$-entropy formula
$(\ref{WW2})$ for the Witten Laplacian on complete Riemannian
manifolds.  Our new observation in Theorem \ref{Thm3} here is that
we can remove the uniform volume lower bound condition
$(\ref{mulow})$ to establish the $W$-entropy formula $(\ref{WW2})$.
This is essentially based on the fact that Theorem \ref{Thm1} and
Theorem \ref{Thm2} hold on complete Riemannian manifolds on which we
need not to assume the uniform volume lower bound condition
$(\ref{mulow})$. Indeed, just after the paper \cite{Li11a} was
accepted, we observed that the uniform volume lower bound condition
$(\ref{mulow})$ can be removed. See the Note added in proof of
\cite{Li11a}.

As an application of Theorem \ref{Thm3}, we prove the following
result, which indicates that the uniform volume lower bound
condition $(\ref{mulow})$  allows us to derive a family of
logarithmic Sobolev inequalities on complete Riemannian manifolds.

\begin{theorem}\label{Thm8a} Let $(M, g, \phi)$ be a weighted complete Riemannian manifold as in Theorem \ref{Thm3}
and satisfying the uniform volume lower bound condition
$(\ref{mulow})$. Then, for any $\tau>0$, the following logarithmic
Sobolev inequality holds: for all $v\in W^{1, 2}(M, \mu)$ with
$\int_M v^2d\mu=1$, we have
\begin{eqnarray*}
\int_M v^2\log v^2 d\mu\leq 4\tau \int_M |\nabla
v|^2-m\left(1+{1\over 2}\log(4\pi\tau)\right)-\mu(\tau),
\end{eqnarray*}
where
\begin{eqnarray*}
\mu(\tau):=\inf\limits\left\{W_m(u, \tau): \
u={e^{-f}\over (4\pi \tau)^{m/2}}, \ \int_M ud\mu=1\right\}
\end{eqnarray*}
is a finite number.
\end{theorem}

\begin{remark}\label{rem1} {\rm Theorem \ref{Thm1} and Theorem \ref{Thm2} are improved versions of Theorem 4.1 and Theorem 4.2 obtained in \cite{Li11a}, where the uniform volume lower bound condition $(\ref{mulow})$ is assumed. When M is a compact Riemannian manifold and $\phi$ is a constant, Theorem \ref{Thm1} and Theorem \ref{Thm2}
are due to Hamilton \cite{H1}, Sheu \cite{Sh}, Hsu \cite{Hsu} and Stroock-Turesky \cite{ST}.
In the case where $L=\Delta$ is the Laplace-Beltrami operator on complete Riemannian manifolds with bounded geometry
condition (in particular, the uniform volume lower bound condition $(\ref{mulow})$ holds for $\mu=\nu$) and for small time $t>0$, Theorem \ref{Thm1} is due to Engoulatov \cite{En}. In \cite{Ko, DG}, the authors proved the gradient estimates in Theorem \ref{Thm1} for $L=\Delta$ on complete Riemannin manifolds with non-negative Ricci curvature.}
\end{remark}

\begin{remark}{\rm In  \cite{LL2}, Songzi Li and the author gave a new proof of the $W$-entropy formula on Riemannian manifolds by using the warped metric product \begin{eqnarray*}
\widetilde{g}=g_M\bigoplus e^{-{2\phi\over m-n}}g_{N}.\label{WPM}
\end{eqnarray*} on $M\times N$, where $(N, g_N)$ is any fixed compact  Riemannian manifold of dimension $m-n$. This gives a natural geometric interpretation of the third term appeared in $(\ref{WW2})$. See Remark $2.2$ in \cite{LL2}. Moreover, we extended the $W$-entropy formula $(\ref{WW2})$ to the heat equation of the time dependent Witten Laplacian on compact Riemannian manifolds on which the metrics $\{g(t), t\in [0, T]\}$ and the potential functions $\{\phi(t), t\in [0, T]\}$ evolve along the $m$-dimensional Perelman  Ricci flow and the conjugate heat equation
\begin{eqnarray*}
{\partial g\over \partial t}&=&-2\left(Ric+\nabla^2\phi-{\nabla\phi\otimes \nabla \phi\over m-n}\right),\\
{\partial \phi\over \partial t}&=&-R-\Delta \phi+{|\nabla \phi|^2\over m-n}.
\end{eqnarray*}
In this case, without assuming $Ric_{m, n}(L)\geq 0$, we prove in \cite{LL2} that the $W$-entropy for the positive solution of the forward heat equation
$$\partial_t u=Lu$$
of the time dependent Witten Laplacian
$L=\Delta-\nabla\phi\cdot\nabla$, is non-increasing in time $t\in [0, T]$.
}
\end{remark}

The rest of this paper is organized as follows.  In Section $2$, we prove the Hamilton type Harnack inequalities (Theorem \ref{Thm0} and Theorem \ref{Thm02}) for the Witten Laplacian on complete Riemanian manifolds.  In Section $3$, we use Theorem \ref{Thm02} to prove the first order gradient estimate of the logarithmic heat kernel of the Witten Laplacian, i.e., Theorem \ref{Thm1}.
In Section $4$, we use a probabilistic approach to prove the gradient estimates of the logarithmic heat kernel of the Witten Laplacian, i.e., Theorem \ref{Thm1} and Theorem \ref{Thm2}. In Section $4$, we prove the $W$-entropy formula and derive the monotonicity and rigidity theorem of the $W$-entropy for the Witten Laplacian on complete Riemannian manifolds, i.e., Theorem \ref{Thm3}. In Section $5$, we prove a family of logarithmic Sobolev inequalities on complete Riemannian manifolds with the uniform volume lower bound condition, i.e., Theorem \ref{Thm8a}.

\section{Hamilton's Harnack inequality}

In this section we give two proofs of Theorem \ref{Thm0} and Theorem \ref{Thm02}.
\medskip

\noindent{\bf Direct Proof of Theorem \ref{Thm0}}. Similarly to \cite{H1}, a direct calculation and the generalized Bochner formula yield
\begin{eqnarray*}
(\partial_t-L){|\nabla u|^2\over u}&=&-{2\over u}|\nabla^2 u-u^{-1}\nabla u\otimes \nabla u|^2-2u^{-1}(Ric+\nabla^2\phi)(\nabla u, \nabla u)\\
&\leq&2Ku^{-1}|\nabla u|^2.
\end{eqnarray*}
Moreover
\begin{eqnarray*}
(\partial_t-L) (u\log (A/u))={|\nabla u|^2\over u}.
\end{eqnarray*}
Let $\psi(t)={1-e^{-2Kt}\over 2K}$. Then
\begin{eqnarray*}
\psi'(t)+2K\psi(t)=1.
\end{eqnarray*} Define
\begin{eqnarray*}
h(x, t):=\psi(t){|\nabla u|^2\over u}-u\log (A/u).
\end{eqnarray*}
We have
\begin{eqnarray*}
(\partial_t-L)h&=&\psi'(t){|\nabla u|^2\over u}+\psi(t)(\partial_t-L){|\nabla u|^2\over u}-(\partial_t-L)u\log (A/u)\\
&\leq& \psi'(t){|\nabla u|^2\over u}+2K\psi(t){|\nabla u|^2\over u}-{|\nabla u|^2\over u}\\
&=& 0.
\end{eqnarray*}
Note that
$$
h(x, 0)=-u\log(A/u)(x, 0)\leq 0.$$

In the case $M$ is compact, the maximum principle yields that $h(x, t)\leq 0$ for all time $t>0$ and $x\in M$. In the case $M$ is a complete non-compact Riemannian manifold with $Ric(L)\geq -K$,  we can give a probabilistic proof to $(\ref{BLH})$ as follows. Indeed, on any complete Riemannian manifold $M$ with $Ric(L)\geq -K$, a previous result due to Bakry \cite{Ba82} says that the $L$-diffusion process on $M$ has infinity lifetime. Let $X_t$ be the $L$-diffusion process on $M$ starting from $X_0=x$. Applying It\^o's formula to $h(X_t, T-t)$, $t\in [0, T]$, we have
\begin{eqnarray*}
h(X_t, T-t)=h(X_0, T)+\int_0^t \nabla h(X_s, T-s)\cdot dW_s+\int_0^t \left(L-{\partial\over \partial t}\right)h(X_s, T-s)ds,
\end{eqnarray*}
where the second term in the right hand side is the It\^o's stochastic integral with respect to a $M$-valued Brownian motion $\{W_s, s\in [0, t]\}$.
In particular, taking $t=T$, we obtain
\begin{eqnarray*}
h(X_T, 0)=h(X_0, T)+\int_0^T \nabla h(X_s, T-s)\cdot dW_s+\int_0^T \left(L-{\partial\over \partial t}\right)h(X_s, T-s)ds.
\end{eqnarray*}
Taking the expectation on both sides, the martingale property of It\^o's integral implies that
\begin{eqnarray*}
E\left[h(X_T, 0)\right]=h(x, T)+E\left[\int_0^T \left(L-{\partial\over \partial t}\right)h(X_s, T-s)ds\right]\geq h(x, T).
\end{eqnarray*}
As $h(y, 0)\leq 0$ for all $y\in M$, we derive that $h(x, T)\leq 0$ for all $T>0$ and $x\in M$.  \hfill $\square$

\medskip

\noindent{\bf Alternative proof of Theorem \ref{Thm0}}. By \cite{Ba, BL}, we know that $Ric(L)\geq -K$ if and only the following version of logarithmic Sobolev inequality holds:  for all $T>0$, $f\in C_b(M)$ with $f>0$,
\begin{eqnarray}
{|\nabla P_T f|^2\over P_T f}\leq {2K\over 1-e^{-2KT}}\left(P_T(f\log f)-P_Tf\log P_Tf\right).\label{RLSI}
\end{eqnarray}
Replacing $f$ by $P_tf$ and using the fact that $0<P_t f(x)\leq A$ for all $x\in M$, one has
\begin{eqnarray*}
{|\nabla P_{T+t} f|^2\over P_{T+t} f}\leq {2K\over 1-e^{-2KT}}\left(P_T(P_t f\log A)-P_{T+t}f\log P_{T+t}f\right).
\end{eqnarray*}
Thus
\begin{eqnarray*}
|\nabla \log P_{T+t} f|^2\leq {2K\over 1-e^{-2KT}}\log (A/P_{T+t}f), \ \ \forall\ t>0.
\end{eqnarray*}
Taking $t\rightarrow 0^+$ we have
\begin{eqnarray*}
|\nabla \log P_{T} f|^2\leq {2K\over 1-e^{-2KT}}\log (A/P_{T}f).
\end{eqnarray*}
This finishes the proof of Theorem \ref{Thm0}. \hfill $\square$

\noindent{\bf Proof of Theorem \ref{Thm02}}.  We can derive the Hamilton Harnack inequality $(\ref{Ham})$ from the Harnack inequality $(\ref{BLH})$ in Theorem \ref{Thm0}. More precisely, using the inequality ${2K\over 1-e^{-2KT}}\leq 2K+{1\over T}$,  we have
\begin{eqnarray*}
|\nabla \log P_{T} f|^2
\leq \left(2K+{1\over T}\right)\log (A/P_{T}f).
\end{eqnarray*}
Note that, if we follow Hamilton \cite{H1} to use the function
$\psi(t)={t\over 2Kt+1}$ (instead of $\psi(t)={1-e^{-2Kt}\over
2K}$), then $\psi'(t)+2K\psi(t)\leq 1$. By the same argument as in
the direct proof of Theorem \ref{Thm0}, we can give a direct proof
of the Harnack inequality $(\ref{Ham})$. \hfill $\square$

\begin{remark}\label{rem3} {\rm As was pointed out in Section 1, Theorem \ref{Thm02} was first proved by Hamilton \cite{H1} for the heat equation of the Laplace-Beltrami operator on compact Riemannian manifolds  and was extended to complete Riemannian manifolds with Ricci curvature bounded from below by Kotschwar \cite{Ko}. For a probabilistic proof of Hamilton's Harnack inequality for $L=\Delta$ on compact or complete Riemannian manifolds, see \cite{AT}. The local version of Theorem \ref{Thm0} for the Laplace-Beltrami operator was proved by Souplet-Zhang \cite{SZ}. In \cite{ATW} (Theorem 5.1),  Arnaudon, Tahlmaier and Wang extended the local estimate of Souplet-Zhang \cite{SZ} to the heat equation of the Laplacian with general drift $L=\Delta+Z$,  where $Z$ is not necessary of gradient form $Z=\nabla \phi$ and the Laplacian comparison theorem for $L=\Delta+Z$ is essentially used. Taking $R\rightarrow \infty$ in \cite{ATW}, we can derive the following estimate for the bounded and positive solution of the heat equation $\partial_t u=Lu$: there exists a constant $C>0$, which depends only on $n={\rm dim}M$, such that for all $x\in M$ and $t>0$, it holds
\begin{eqnarray}
|\nabla \log u|^2\leq C\left({1\over t}+K\right)\left[1+\log(A/u)\right].\label{ATMW}
\end{eqnarray}
where $A:=\sup\limits\{u(t, x): x\in M, t\geq 0\}$.
}
\end{remark}

As a corollary of the above theorem, we can derive the following gradient estimate and the Liouville theorem due to Brighton \cite{Br} for bounded $L$-harmonic functions.

\begin{corollary}\label{BHT} (i.e., Corollary \ref{cor1})) Suppose that $Ric(L)=Ric+\nabla^2\phi\geq -K$ for some constant $K\geq 0$. Then every bounded  $L$-harmonic function $u$ satisfies
\begin{eqnarray}
|\nabla \log u|\leq 2K\log(A/(u-\inf\limits u)), \label{grad-1}
\end{eqnarray}
where $A=\sup\limits u-\inf\limits u$. In particular, if $Ric(L)=Ric+\nabla^2\phi\geq 0$, then every bounded $L$-harmonic function must be constant.
\end{corollary}
{\it Proof}. Applying $(\ref{Ham})$ to $u(t, x)=u(x)-\inf\limits u$, $x\in M$, $t>0$, we have
\begin{eqnarray*}
|\nabla\log u(x)|^2\leq \left(2K+{1\over T}\right)\log\left({A/(u(x)-\inf\limits u)}\right),\ \ \forall T>0.
\end{eqnarray*}
Taking $T\rightarrow \infty$, we derive $(\ref{grad-1})$. In particular, if $K=0$, we have $|\nabla \log u|=0$. This finishes the proof of Corollary \ref{BHT}. \hfill $\square$

\begin{remark}\label{rem4} {\rm In the case $L=\Delta$, Yau \cite{Yau} proved that, on complete Riemannian manifolds with Ricci curvature bounded below by $-K$, i.e., $Ric\geq K$, every harmonic function which is bounded below by a constant satisfies the gradient estimate
\begin{eqnarray*}
|\nabla \log u|\leq \sqrt{(n-1)K}(u-\inf\limits u).\label{grad-2}
\end{eqnarray*}
In particular, Yau \cite{Yau} proved the Strong Liouville theorem holds: On complete Riemannian manifolds with non-negative Ricci curvature, every positive (and hence bounded) harmonic function must be constant. See also \cite{SY}. In \cite{Li05}, the author extended Yau's gradient estimate and Strong Liouville theorem to the $L$-harmonic functions on complete Riemnnian manifolds via the Bakry-Emery Ricci curvature: Suppose that $Ric_{m, n}(L)\geq -K$ for some constant $K\geq 0$. Then every $L$-harmonic function $u$, which is bounded from below, satisfies the following gradient estimate
\begin{eqnarray*}
|\nabla \log u|\leq \sqrt{(m-1)K}(u-\inf\limits u).\label{grad-3}
\end{eqnarray*}
In particular, if $Ric_{m, n}(L)\geq 0$, then all positive (and hence bounded) $L$-harmonic functions must be constant. In \cite{Li07}, the local gradient estimates was proved for $L$-harmonic functions: If $Ric_{m, n}(L)\geq -K$, then every positive $L$-harmonic function $u$ on $M$ satisfies: for all $o\in M$, and $R>0$, we have
\begin{eqnarray*}
|\nabla \log u(x)|\leq \sqrt{(m-1)K}u(x)+ C(K, m)(R^{-1/2}+R^{-3/2}), \ \ \forall x\in B(o, R).\label{grad-4}
\end{eqnarray*}
In \cite{Br}, Brighton proved that, on complete Riemannian manifolds with $Ric+\nabla^2\phi\geq 0$, every bounded $L$-harmonic function must be constant, i.e., Corollary \ref{cor1}. Our method for the proof of Corollary \ref{cor1} is different from \cite{Br}.
For more recent results, see \cite{MW, MW2}.}
\end{remark}
\bigskip

\section{Proof of Theorem \ref{Thm1}}

In this section, we use Hamilton's Harnack inequality to  prove Theorem \ref{Thm1}. This proof is similar to the one given in \cite{Ko, DG} for the gradient estimate of the logarithmic heat kernel of the Laplace-Beltrami operator on complete Riemannian manifolds with non-negative Ricci curvature.

To prove Theorem \ref{Thm1}, we need the following estimates on the heat kernel of the Witten Laplacian on complete Riemannian manifolds with weighted volume measure.

\begin{proposition}\label{prop3} Suppose that there exist some constants $m\geq n$, $m\in \mathbb{N}$ and $K\geq 0$ such that $Ric_{m, n}(L)\geq -K$. Then, for any small $\varepsilon>0$, there exist some constants $C_i=C_i(m, n, K, \varepsilon)>0$, $i=1, 2$, such that for all $x, y\in M$ and $t>0$,
\begin{eqnarray}
p_t(x, y)&\leq& {C_1\over \mu(B_y(\sqrt{t}))}\exp\left(-{d^2(x, y)\over 4(1+\varepsilon)t}+\alpha \varepsilon Kt\right)\nonumber\\
& &\hskip1cm \times \left({d(x, y)+\sqrt{t}\over \sqrt{t}}\right)^{m/2} \exp\left({\sqrt{(m-1)K}d(x, y)\over 2}\right),\label{heat-upp}
\end{eqnarray}
where $\alpha$ is a constant depending only on $m$, and
\begin{eqnarray}
p_t(x, y)\geq C_2e^{-(1+\varepsilon)\lambda_{K,
m}t}\mu^{-1}(B_y(\sqrt{t}))\exp\left(-{d^2(x, y)\over
4(1-\varepsilon)t}\right)\left[{\sqrt{K}d(x, y)\over \sinh \sqrt{K}d(x,
y)}\right]^{m-1\over 2},\label{heat-lower}
\end{eqnarray}
where
$$\lambda_{K, m}={(m-1)^2K\over 8}.$$
\end{proposition}
{\it Proof}. The lower bound estimate $(\ref{heat-lower})$ has proved in \cite{Li11a}. By Theorem 5.4 in \cite{Li05}, for all $\varepsilon>0$, there exists a constant $C_3=C_3(m, n, K, \varepsilon)>0$
\begin{eqnarray}
p_t(x, y)\leq {C_3\over \sqrt{\mu(B_x(\sqrt{t}))\mu(B_y(\sqrt{t}))}}\exp\left(-{d^2(x, y)\over 4(1+\varepsilon)t}+\alpha \varepsilon K t\right).\label{heat-upp-2}
\end{eqnarray}
Similarly to \cite{Li05} (Step 2 of p. 1324), using the Bishop-Gromov relative volume comparison theorem \cite{Q, Li05}, we have
\begin{eqnarray*}
\mu(B_y(\sqrt{t}))&\leq& \mu\left(B_x(d(x, y)+\sqrt{t})\setminus B_x(d(x, y)-\sqrt{t})\right)\\
&\leq& \mu(B_x(\sqrt{t})){V_{m, K}(d(x, y)+\sqrt{t})-V_{m, K}(d(x, y)-\sqrt{t})\over V_{m, K}(\sqrt{t})}\\
&\leq&\mu(B_x(\sqrt{t})){V_{m, K}(d(x, y)+\sqrt{t})\over V_{m, K}(\sqrt{t})},
\end{eqnarray*}
where $V_{m, K}(r)$ denotes the volume of geodesic balls of radius $r$ in the $m$-dimensional hyperbolic space form $\mathbb{H}^m(K)$ of constant sectional curvature $K/m-1$.
Using again the Bishop-Gromov volume comparison theorem on $\mathbb{H}^m(K)$, we have
\begin{eqnarray}
{\mu(B_y(\sqrt{t}))\over \mu(B_x(\sqrt{t}))}
&\leq&{V_{m, K}(d(x, y)+\sqrt{t})\over V_{m, K}(\sqrt{t})}\nonumber\\
&\leq& \left({d(x, y)+\sqrt{t}\over \sqrt{t}}\right)^m \exp[\sqrt{(m-1)K}d(x, y)].\label{vol-mu1}
\end{eqnarray}
By $(\ref{heat-upp-2})$ and $(\ref{vol-mu1})$, we can obtain the upper bound estimate $(\ref{heat-upp})$. \hfill $\square$

\medskip

\noindent{\bf Proof of Theorem \ref{Thm1}}.
Fix $T>0$,  and let $u(t, x)$ be a positive and bounded solution to the heat equation
$\partial_t u=L u$, $t\in (0, t_1)$. Let
$$
A:=\sup\limits\{u(t, x): 0\leq t\leq t_1, x\in M\}.$$ By
Hamilton's Harnack inequality $(\ref{Ham})$, we have
\begin{eqnarray}
t|\nabla \log u(t, x)|^2\leq (1+2Kt)\log(A/u(t, x)), \ \ \ \forall
(t, x)\in [0, t_1]\times M. \label{EE1}
\end{eqnarray}
Let $s\in (0, T]$, $y\in M$, $t_1=s/2$ and $u(t, x)=p_{s/2+t}(x,
y)$. By $(\ref{EE1})$, $(\ref{heat-upp})$ and $(\ref{heat-lower})$,  we have
\begin{eqnarray*}
& &t|\nabla_x \log p_{s/2+t}(x, y)|^2\leq C_{K, m, T}(1+Kt) \\
& &\hskip1cm \times \left[1+{d(x, y)\over \sqrt{t}}+\log\left({C_1\over
C_2}{\mu(B(y, \sqrt{s/2+t}))\over \mu(B(y, \sqrt{s/2}))} \exp\left(C_3{d^2(x,
y)\over s/2+t}+C_4d(x, y)\right)\right)\right].
\end{eqnarray*}
In particular, taking $t=s/2$ and changing $s$ by $t$ we get
\begin{eqnarray*}
& &{t\over 2} |\nabla_x \log p_{t}(x, y)|^2\leq  C_{K, m, T}(1+K{t\over
2})\\
& &\hskip1cm \times \left[1+{d(x, y)\over \sqrt{t}}+\log\left({C_1\over C_2}{\mu(B(y, \sqrt{t}))\over \mu(B(y, \sqrt{t/2}))}
\exp\left(C_3{d^2(x, y)\over t}+C_4d(x, y)\right)\right)\right].
\end{eqnarray*}
By the Bishop-Gromov volume comparison $(\ref{E6})$, we derive
\begin{eqnarray*}
t|\nabla_x \log p_{t}(x, y)|^2\leq C_{K, m, T}\left(1+{d^2(x,
y)\over t}+{d(x, y)\over \sqrt{t}}+d(x, y)\right),
\end{eqnarray*}
which yields
\begin{eqnarray*}
|\nabla_x \log p_{t}(x, y)|\leq C_{K, m, T}\left({d(x, y)\over
t}+{1\over \sqrt{t}}\right).
\end{eqnarray*}
This finishes the analytic proof of Theorem \ref{Thm1}. \hfill $\square$

\bigskip

\section{Gradient estimates of the heat kernel}

In this section we use a probabilistic approach to prove Theorem \ref{Thm1} and Theorem \ref{Thm2}. Our proof is inspired  by previous works initialed by Sheu \cite{Sh}, and developed by Hsu \cite{Hsu} and Engoulatov \cite{En}.

Let $O(M)$ be the orthogonal frame bundle over $M$, $\pi: O(M)\rightarrow M$ the canonical projection map. Let $H_1, \ldots, H_n$ be the canonical horizonal vector fields on $O(M)$. Let $B_t$ be the standard Brownian motion on $\mathbb{R}^n$. Following Malliavin \cite{Ma74, Ma97}, we define the horizontal $L$-diffusion process $U_t$ on $O(M)$ by the following Stratonovich SDE on $O(M)$:
\begin{eqnarray}
dU_t=\sum\limits_{i=1}^n H_i(U_t)\circ dB^i_t-\nabla^H \phi(U_t)dt,\label{U}
\end{eqnarray}
where $\nabla^H \phi$ denotes the horizontal lift of the gradient vector field $\nabla \phi$ on $O(M)$, which is the unique horizontal vector field on $O(M)$ such that
\begin{eqnarray*}
\pi_*(\nabla^H \phi)=\nabla \phi.
\end{eqnarray*}
Let
\begin{eqnarray*}
X_t=\pi(U_t),  \ \ \forall t>0.
\end{eqnarray*}
Then $X_t$ is a diffusion process on $M$ with infinitesimal generator $L$. Moreover, we have
\begin{eqnarray}
dX_t=U_t\circ dB_t-\nabla \phi(X_t)dt.\label{X}
\end{eqnarray}

Let  $\Delta_{O(M)}=\sum\limits_{i=1}^n H_i^2$ be the horizontal Laplace-Beltrami on $O(M)$. The horizontal Witten Laplacian is deined by
$$L_{O(M)}=\Delta_{O(M)}-\nabla^H \phi\cdot\nabla^H.$$

Fix $T>0$. Let $\mathbb{P}_{x, y, T}$ be the conditional $L$-diffusion processes starting from $x$ and ending at $y$ t time $T$, i.e.,
\begin{eqnarray*}
\mathbb{P}_{x, y, T}(\cdot)=\mathbb{P}_{x}(\cdot|X_T=y).
\end{eqnarray*}
Let
\begin{eqnarray*}
J(t, u)=\log p_{T-t}(\pi u, y), \ \ \ \forall u\in O(M),  t\in [0, T].
\end{eqnarray*}
By the Grirsanov transformation, we have
\begin{eqnarray*}
\left. {d\mathbb{P}_{x, y, T}\over d\mathbb{P}_{x}}\right|_{\mathcal{F}_t}={p_{T-t}(X_t, y)\over p_T(x, y)}=\exp\{J(t, U_t)-J(0, U_0)\},\ \ \ \forall t<T.
\end{eqnarray*}
where $\mathcal{F}_t=\sigma(X_s, s\leq t)$ is the natural $\sigma$-filed generated by $X_t$. By It\^o's formula, we have
\begin{eqnarray*}
dJ(t, U_t)={\partial J\over \partial t}(t, U_t)dt+\sqrt{2}\langle \nabla^H J(t, U_t), dB_t\rangle+L_{O(M)}J(t, U_t)dt,
\end{eqnarray*}
Note that $J(t, u)$ satisfies the Hamilton-Jacobi equation
\begin{eqnarray}
{\partial J\over \partial t}+L_{O(M)}J+|\nabla^H J|^2=0,\label{HJJ}
\end{eqnarray}
which yields
\begin{eqnarray*}
\left. {d\mathbb{P}_{x, y, T}\over d\mathbb{P}_{x}}\right|_{\mathcal{F}_t}=\exp\left[\int_0^t \langle\sqrt{2} \nabla^H J(s, U_s), dB_s\rangle-\int_0^T |\nabla^H J(s, U_s)|^2ds\right],\ \ \ \forall t<T.
\end{eqnarray*}
It follows, by Girsanov's theorem, that the process
\begin{eqnarray}
b_t=B_t-\sqrt{2}\int_0^t \nabla^H J(s, U_s)ds, \ \ \ s\in [0, T),\label{b}
\end{eqnarray}
is a Brownian motion under $\mathbb{P}_{x, y, T}$. Substituting $(\ref{b})$ into $\ref{U})$ and $(\ref{X})$, we have the following lemma.

\begin{lemma}\label{lem1} The conditional $L$-diffusion process $X_t$ on $M$ staring from $x$ with terminal end $y$ at the time $T$, and its horizontal lift $U_t$, are the solutions of the following Stratonovich SDEs:
\begin{eqnarray*}
dU_t&=&\sum\limits_{i=1}^n H_i(U_t)\circ\left[\sqrt{2}db_t^i+2H_i\log p_{T-t}(U_t, y)dt\right]-\nabla^H \phi(U_t)dt,\\
dX_t&=&U_t\circ \left[\sqrt{2}db_t+2\nabla^H \log p_{T-t}(U_t, y)dt\right]-\nabla\phi(X_t)dt,
\end{eqnarray*}
where $b_t$ is a standard Brownian motion on $\mathbb{R}^n$ under $\mathbb{P}_{x, y, T}$.
\end{lemma}

\begin{proposition}\label{prop1} For all $0<t<T$, we have
\begin{eqnarray}
\mathbb{E}[d|\nabla^H J(t, U_t)|^2]=2\mathbb{E}\left[|\nabla^H\nabla^H J|^2+(Ric+\nabla^2\phi)(\nabla^H J, \nabla^H J)\right].\label{E1}
\end{eqnarray}
\end{proposition}
{\it Proof}. By Lemma \ref{lem1} and It\^o's formula, we have
\begin{eqnarray*}
d|\nabla^H J(t, U_t)|^2=\langle \nabla^H |\nabla^H J|^2, \sqrt{2}db_t+2\nabla^H Jdt\rangle+(\partial_t+L_{O(M)})|\nabla^H J|^2dt.
\end{eqnarray*}
By the generalized Bochner-Weitzenb\"ock formula, the action of the horizontal Witten Laplacian $L_{O(M)}$ on $|\nabla^H J|^2$ is given by
\begin{eqnarray*}
L_{O(M)}|\nabla^H J|^2=2\langle \nabla^H L_{O(M)}J, \nabla^H J\rangle+2|\nabla^H\nabla^H J|^2+2(Ric+\nabla^2\phi)(\nabla^H J, \nabla^H J).
\end{eqnarray*}
On the other hand, using the Hamilton-Jacobi equation $(\ref{HJJ})$, we have
\begin{eqnarray*}
\partial_t |\nabla^H J|^2=2\langle \partial_t\nabla^H J, \nabla^H J\rangle=-2\langle L_{O(M)}\nabla^H J, \nabla^H J\rangle-2\langle\nabla^H|\nabla^H J|^2, \nabla^H J\rangle.
\end{eqnarray*}
Combining this with the previous calculation, we have
\begin{eqnarray*}
d|\nabla^H J(t, U_t)|^2=\langle \nabla^H |\nabla^H J|^2, \sqrt{2}db_t\rangle+2\left[|\nabla^H\nabla^H J|^2+(Ric+\nabla^2\phi)(\nabla^H J, \nabla^H J)\right]dt.
\end{eqnarray*}
Taking expectation with respect to $\mathbb{E}_{x, y, T}$, we complete the proof of Proposition \ref{prop1}. \hfill $\square$
\begin{proposition}\label{prop2} Suppose that $Ric+\nabla^2\phi\geq -K$. Then for all $T>0$, we have
\begin{eqnarray}
|\nabla\log p_{T}(x, y)|^2\leq 2\left({1\over T}+K\right)\mathbb{E}\left[\log {p_{T/2}(X_{T/2}, y)\over p_T(x, y)}\right].\label{E3}
\end{eqnarray}
\end{proposition}
{\it Proof}.  Based on Proposition \ref{prop1}, the proof of Proposition \ref{prop2} is similar to the argument used in \cite{Hsu, En}. For the convenience of the reader, we give the detail here. Under the condition $Ric+\nabla^2\phi\geq -K$, integrating $(\ref{E1})$ from $0$ to $t$, we have
\begin{eqnarray*}
\mathbb{E}[|\nabla^HJ(t, U_t)|^2]-|\nabla^H J(0, U_0)|^2\geq -2KE\left[\int_0^t |\nabla^H J(s, U_s)|^2ds\right].
\end{eqnarray*}
Integrating again with respect to $t$ we can obtain
\begin{eqnarray}
{T\over 2}|\nabla^H J(0, U(0)|^2\leq (1+KT)\mathbb{E}\left[\int_0^{T/2}|\nabla^H J(s, U_s)|^2ds\right].\label{E2}
\end{eqnarray}
By It\^o's formula and the Hamilton-Jacobi equation $(\ref{HJJ})$, we have
\begin{eqnarray*}
dJ(t, U_t)&=&(\partial_t+L_{O(M)})Jdt+\langle \sqrt{2}\nabla^H J, db_t+\sqrt{2}\nabla^H Jdt\rangle\\
&=&\langle \sqrt{2}\nabla^H J, db_t\rangle+|\nabla^H J|^2dt.
\end{eqnarray*}
Taking expectation with respect to $\mathbb{E}_{x, y, T}$ and integrating from $t=0$ to $t=T/2$, we get
\begin{eqnarray}
\mathbb{E}[J(T/2, U_{T/2})]-J(0, U_0)=\mathbb{E}\left[\int_0^{T/2}|\nabla^H J|^2dt\right].\label{E4}
\end{eqnarray}
Combining $(\ref{E2})$ with $(\ref{E4})$, we finish the proof of $(\ref{E3})$. \hfill $\square$

\bigskip

\noindent{\bf Probabilistic proof of Theorem \ref{Thm1}}. By Proposition \ref{prop3}, we have
\begin{eqnarray}
{p_{t/2}(X_{t/2}, y)\over p_{t}(x, y)}&\leq& C_1(K, m, T) {\mu(B_y(\sqrt{t}))\over \mu(B_y(\sqrt{t/2}))} \exp\left({d^2(x, y)\over
4(1-\varepsilon)t}+{\sqrt{(m-1)K}d(X_{t/2}, y)\over 2}\right)\nonumber\\
& &\hskip2cm \times  \left({d(X_{t/2}, y)+\sqrt{t}\over \sqrt{t}}\right)^{m/2}\left[{\sqrt{K}d(x, y)\over \sinh \sqrt{K}d(x,
y)}\right]^{-{m-1\over 2}}.\label{E5}
\end{eqnarray}
The Bishop-Gromov volume comparison theorem (\cite{Q, Li05}) implies
\begin{eqnarray}
{\mu(B_y(\sqrt{t}))\over \mu(B_y(\sqrt{t/2}))}\leq 2^{m/2}
\exp(\sqrt{m-1}Kt)\leq 2^{m/2}\exp(\sqrt{m-1}KT).\label{E6}
\end{eqnarray}
Using  $(\ref{E5})$, $(\ref{E6})$ and the elementary inequality $x/\sinh x\geq e^{-x}$, we obtain
\begin{eqnarray*}
\log\left[{p_{t/2}(X_{t/2}, y)\over p_{t}(x, y)}\right]
&\leq&C_{2}(K, m, T)\left[1+{d(X_{t/2}, y)\over \sqrt{t}}+{d^2(x, y)\over t}+d(X_{t/2}, y)+d(x, y)\right]\nonumber\\
&\leq&C_{3}(K, m, T)\left[1+{d^2(X_{t/2}, y)\over t}+{d^2(x, y)\over t}+d(X_{t/2}, y)+d(x, y)\right].\label{AE1}
\end{eqnarray*}
By It\^o's formula, we can prove that (cf. \cite{Li11a})
\begin{eqnarray*}
\mathbb{E}_{x, y, T}[d^k(x_{t/2}, y)]\leq (d^k(x, y)+1)e^{C_4(K, m)t}, \ \ \ \forall t>0, \ k=1, 2.
\end{eqnarray*}
Taking expectation with respect to $\mathbb{E}_{x, y, T}$ in $(\ref{AE1})$, we have
\begin{eqnarray*}
\mathbb{E}_{x, y, T}\left[\log {p_{t/2}(X_{t/2}, y)\over p_{t}(x, y)}\right]\leq
C_{5}(K, m, T)\left[1+{d^2(x, y)+1\over t}+d(x, y)\right].
\end{eqnarray*}
By Proposition \ref{prop2}, we have
\begin{eqnarray*}
|\nabla \log p_t(x, y)|^2&\leq& C_{6}(K, m, T)\left[{1\over t}+{d^2(x, y)+1\over t^2}+{d(x, y)\over t}\right]\\
&\leq& C_7(K, m, T)\left[{d(x, y)\over t}+{1\over \sqrt{t}}\right]^2.
\end{eqnarray*}
The proof of Theorem \ref{Thm1} is completed. \hfill $\square$

\bigskip

\noindent{\bf Proof of Theorem \ref{Thm2}}. The proof is similarly to Sheu \cite{Sh} and Hsu \cite{Hsu}.
Let $I=\{i_1, \ldots, i_N\}$ ,and denote $H_I=H_{i_1}\ldots H_{i_N}$, then by Lemma \ref{lem1} and repeatedly using It\^o's formula, we have

\begin{eqnarray}
dH_I J=\langle \nabla^H H_I J, \sqrt{2}db_t\rangle+\left[\partial_t H_I J+L_{O(M)}H_I J+2\langle \nabla^H H_I J, \nabla^H J\rangle\right]dt.\label{E7}
\end{eqnarray}
Substituting the Hamilton-Jacobi equation $(\ref{HJJ})$ into $(\ref{E7})$, we have
\begin{eqnarray}
dH_IJ =\langle \nabla^H H_I J, \sqrt{2}db_t\rangle+(F_I+G_I)dt.\label{E8}
\end{eqnarray}
where
\begin{eqnarray*}
F_I=[L_{O(M)}, H_I]J=[\Delta_{O(M)}, H_I]J-[\nabla^H \phi, H_I]J,
\end{eqnarray*}
and
\begin{eqnarray*}
G_I=2\langle \nabla^H H_I J, \nabla^H J\rangle-H_I|\nabla^H J|^2.
\end{eqnarray*}
By $(\ref{E4})$ and  $(\ref{E5})$ , we have
\begin{eqnarray}
\mathbb{E}\left[\int_0^{T/2}|\nabla^H J(t, U_t)|^2dt\right]\leq CT\left[{d(x, y)\over T}+{1\over \sqrt{T}}\right]^2, \label{E9}
\end{eqnarray}
and by Theorem \ref{Thm1}, we have
\begin{eqnarray}
|H_iJ(0, U_0)|\leq C\left[{d(x, y)\over T}+{1\over \sqrt{T}}\right].\label{E10}
\end{eqnarray}
Note that
\begin{eqnarray*}
[\nabla^H\phi, H_I]J=\sum\limits_{k=1}^n[(H_k \phi)H_k, H_I]J=\sum\limits_{k=1}^n (H_k \phi[H_k, H_I]J-H_IH_k\phi H_kJ),
\end{eqnarray*}
and by Cartan's structure theorem, we have
\begin{eqnarray*}
[H_i, H_j]=V_{ij}, \ [H_i, V_{jk}]=\Omega_{jk}^{il}H_l, \ [V_{ij}, V_{kl}]=c_{ij, kl}^{ab}V_{ab},
\end{eqnarray*}
where $V_{ij}$ are the vertical fields on $O(M)$, $\Omega$ is the $o(n)$-valued curvature form on $O(M)$, and $c_{ij, kl}^{ab}$ are the structure constants of $o(n)$. The assumption of Theorem \ref{Thm2} says that the $L^\infty$-norms $\|\nabla^k\Omega\|_{\infty}$ are bounded for all $k=0, 1, \ldots, N$, and the $L^\infty$-norms $\|\nabla^{k}\phi\|_{\infty}$ are bounded for $k=1, \ldots, N+1$. This yields
\begin{eqnarray}
|F_I|&\leq& C_1|\nabla^H J|+\ldots +C_1|(\nabla^H)^{|I|}J|,\label{F}\\
|G_I|&\leq& C_1|\nabla^H J|+\ldots +C_1|(\nabla^H)^{|I|}J|+C_1|(\nabla^H)^{|I|}J|^2, \label{G}
\end{eqnarray}
 where $C_1>0$ is a constant depending only on $\|\nabla^i\Omega\|_\infty$ and  $\|\nabla^{i}\phi\|_\infty$, $i=1, \ldots, |I|$. By induction, and using $(\ref{E8})$, $(\ref{E9})$, $(\ref{E10})$, $(\ref{F})$ and $(\ref{G})$, we can prove that for all $I$ with $2\leq |I|\leq N$,
\begin{eqnarray}
\mathbb{E}\left[\int_0^{T/2}|H_I J(t, U_t)|^2dt\right]\leq C_2\left[{d(x, y)\over T}+{1\over \sqrt{T}}\right]^{2(|I|-1)}, \label{E11}
\end{eqnarray}
and for all $I$ with $|I|\leq N$,
\begin{eqnarray}
|H_I J(0, U_0)|\leq C_3\left[{d(x, y)\over T}+{1\over \sqrt{T}}\right]^{|I|}. \label{E12}
\end{eqnarray}
where $C_2$ and $C_3$ are two constants depending on $K, m, n, T$ and $C_1$. The proof of Theorem \ref{Thm2} is completed. For the detail of the above induction argument for $L=\Delta$ on compact Riemannian manifolds, see Hsu \cite{Hsu}. \hfill $\square$

\section{Proofs of Theorem \ref{Thm3}}

Let \begin{eqnarray*}
H(u)=-\int_M u\log u d\mu.\label{Ent-0}
\end{eqnarray*}
be the Boltzmann-Shannon-Nash entropy of the heat equation of the Witten Laplacian on $(M, g, \mu)$.

To prove Theorem \ref{Thm3}, we need the
following entropy dissipation formulas for the Witten Laplacian on
complete Riemannian manifolds.

\begin{theorem}(\cite{Li11a})\label{Thm6a} Let $M$ be a complete Riemannian manifold, and $\phi\in C^2(M)$. Suppose that there exists a constant $K\geq 0$ such that
\begin{eqnarray*}
Ric(L)\geq -K.
\end{eqnarray*}
Then, for any positive solution to the heat equation
\begin{eqnarray*}
\partial_t u=L u
\end{eqnarray*}
with $u(\cdot, 0)\in L^1(\mu)\cap L^4(d(\cdot, o)d\mu)$, where $o\in M$ is any fixed point, we have
\begin{eqnarray}
{d\over d t} H(u)=-\int_M Lu\log u d\mu.\label{Ent-1}
\end{eqnarray}
\end{theorem}

\begin{theorem}\label{Thm7} Let $(M, g)$ be an $n$-dimensional complete Riemannian manifold. Suppose that $\|\nabla^k{\rm Riem}\|\leq C_k$ for some constant $C_k>0$, $0\leq k\leq 3$,
$\phi\in C^{4}(M)$ with $\nabla\phi \in C_b^{3}(M)$, and there exists some $m\geq n$, $K\geq 0$ such that \ $Ric_{m, n}(L)\geq -K$.
Then, for the fundamental solution of the heat equation $\partial_t u=Lu$, we have
\begin{eqnarray}
{d^2\over d t^2} H(u)=-\int_M \left({|Lu|^2\over u}-\langle \nabla Lu, {\nabla u\over u}\rangle\right) d\mu,\label{Ent-2}
\end{eqnarray}
and
\begin{eqnarray}
{d^2\over d t^2} H(u)=-2\int_M (|\nabla^2\log u|^2+Ric(L)(\nabla
\log u, \nabla \log u))ud\mu.\label{Ent-3}
\end{eqnarray}
\end{theorem}
{\it Proof}. Based on Theorem \ref{Thm6a}, we can prove $(\ref{Ent-2})$ by the same argument as used in the proof of Theorem 4.3 in \cite{Li11a}.
Integrating by part yields
\begin{eqnarray*}
\int_M {|Lu|^2\over u}d\mu&=&-\int_M \left\langle \nabla\left({Lu\over u}\right), {\nabla u\over u}\right\rangle ud\mu\\
&=&-\int_M \langle \nabla Lu, \nabla\log u\rangle d\mu+\int_M Lu |\nabla \log u|^2d\mu.
\end{eqnarray*}
Hence
\begin{eqnarray*}
{d^2\over dt^2}H(u)&=&2\int_M \left\langle \nabla Lu, {\nabla u\over u}\right\rangle d\mu-\int_M Lu |\nabla \log u|^2 d\mu\\
&=&2\int_M \left\langle \nabla Lu, {\nabla u\over u}\right\rangle d\mu-\int_M u L|\nabla \log u|^2 d\mu.
\end{eqnarray*}
By the Bochner formula for the Witten Laplacian, we have (see \cite{Li11a})
\begin{eqnarray*}
L|\nabla \log u|^2=2\langle\nabla L \log u, \nabla \log
u\rangle+2|\nabla^2\log u|^2+2Ric(L)(\nabla \log u, \nabla \log u)
\end{eqnarray*}
Therefore
\begin{eqnarray*}
{d^2\over dt^2}H(u)&=&2\int_M \left\langle \nabla Lu, {\nabla u\over u}\right\rangle d\mu
-2\int_M \langle \nabla L \log u, \nabla u\rangle d\mu\\
& &\ \ -2\int_M (|\nabla^2\log u|^2+Ric(L)(\nabla
\log u, \nabla \log u))ud\mu.
\end{eqnarray*}
Again, integrating by parts yields
\begin{eqnarray*}
\int_M <\nabla L \log u, \nabla u>d\mu=-\int_M L\log u\cdot Lu d\mu=\int_M \langle\nabla \log u, \nabla Lu\rangle d\mu.
\end{eqnarray*}
Thus
\begin{eqnarray*}
{d^2\over dt^2}H(u)=-2\int_M (|\nabla^2\log u|^2+Ric(L)(\nabla
\log u, \nabla \log u))ud\mu.
\end{eqnarray*}
The proof of Theorem \ref{Thm7} is completed. \hfill $\square$

\begin{remark}\label{rem6} {\rm  Theorem \ref{Thm7} is an improvement of Theorem 4.3 in \cite{Li11a}. The entropy dissipation formulas $(\ref{Ent-1})$ and $(\ref{Ent-2})$ play an important role in the study of the convergence of the heat kernel measure $p_t(\cdot, y)d\mu(y)$ towards the equilibrium measure $\mu$. In the case where $M$ is compact and $\phi\in C^3(M)$, it is well-known that the entropy dissipation formulas $(\ref{Ent-1})$ and $(\ref{Ent-2})$ hold. However, as far as we know, even in the non-weighted case, it seems that there is no explicit result in the literature (except our previous paper \cite{Li11a}) which ensures the entropy dissipation formulas $(\ref{Ent-1})$ and $(\ref{Ent-2})$ on complete non-compact Riemannian manifolds. Theorem \ref{Thm6a} and Theorem \ref{Thm7} give us natural geometric conditions on $(M, g, \phi)$ to ensure  $(\ref{Ent-1})$ and $(\ref{Ent-2})$. }
\end{remark}

\medskip

\noindent{\bf Proof of Theorem \ref{Thm3}}. We simplify the  proof of Theorem 2.3 in \cite{Li11a} as follows. Let
\begin{eqnarray*}
H_m(u, t)=-\int_M u\log u d\mu-{m\over 2}\left(1+\log(4\pi t)\right).
\end{eqnarray*}
By Theorem \ref{Thm6a} and Theorem \ref{Thm7}, we have
\begin{eqnarray}
{d\over dt}H_m(u, t)=-\int_M u L\log u d\mu-{m\over 2t}=\int_M \left(|\nabla \log u|^2-{m\over 2t}\right)ud\mu,\label{Ent-4}
\end{eqnarray}
and
\begin{eqnarray}
{d^2\over dt^2}H_m(u, t)=-2\int_M (|\nabla^2\log u|^2+Ric(L)(\nabla
\log u, \nabla \log u))ud\mu+{m\over 2t^2}.\label{Ent-5}
\end{eqnarray}
By the definition formula of the $W$-entropy, we have
\begin{eqnarray}
W_m(u, t)={d\over dt}(tH_m(u, t))=H_m(u, t)+t{d\over dt}H_m(u, t),\label{Ent-6}
\end{eqnarray}
and
\begin{eqnarray}
{d\over dt}W_m(u, t)=2{d\over dt}H_m(u, t)+t{d^2\over dt^2}H_m(u, t).\label{Ent-7}
\end{eqnarray}
Substituting $(\ref{Ent-0})$ and $(\ref{Ent-4})$ into $(\ref{Ent-6})$, and using $u={e^{-f}\over (4\pi t)^{m/2}}$, we obtain
\begin{eqnarray*}
W_m(u, t)&=&-\int_M u\log u d\mu-{m\over 2}\left(1+\log(4\pi t)\right)-t\int_M uL\log u d\mu-{m\over 2}\\
&=&\int_M \left[t|\nabla \log u|^2-\log u-m-{m\over 2}\log (4\pi t)\right]ud\mu\\
&=&\int_M \left[t|\nabla f|^2+f-m\right]{e^{-f}\over (4\pi t)^{m/2}}d\mu.
\end{eqnarray*}
This proves $(\ref{WW1})$ in Theorem \ref{Thm3}. Substituting $(\ref{Ent-4})$ and $(\ref{Ent-5})$ into $(\ref{Ent-7})$, we obtain
\begin{eqnarray*}
{d\over dt}W_m(u, t)=-2\int_M u L\log u d\mu-2\int_M t(|\nabla^2\log u|^2+Ric(L)(\nabla
\log u, \nabla \log u))ud\mu-{m\over 2t}.\label{W-3}
\end{eqnarray*}
Noting that
\begin{eqnarray*}
2t|\nabla^2\log u|^2+{m\over 2t}=2t\left|\nabla^2 \log u-{g\over 2t}\right|^2+{m-n\over 2t}+2\Delta \log u,
\end{eqnarray*}
and by the integration by part formula
\begin{eqnarray*}
\int_M \Delta \log u ud\mu&=&\int_M (L\log u+\nabla \phi \cdot \nabla \log u) u d\mu\\
&=&-\int_M |\nabla \log u|^2 ud\mu+\int_M \nabla \phi\cdot \nabla \log u u d\mu,
\end{eqnarray*}
we have
\begin{eqnarray*}
{d\over dt}W_m(u, t)&=&-2t\int_M \left|\nabla^2\log u-{g\over 2t}\right|^2ud\mu+2\int_M \nabla \phi\cdot \nabla\log u  ud\mu\\
& &-2\int_M tRic(L)(\nabla \log u, \nabla \log u)ud\mu-{m-n\over 2t}.
\end{eqnarray*}
By the same argument as in \cite{Li11a}, the above identity yields $(\ref{WW2})$. The monotonicity of $W$ follows and the rigidity part of Theorem \ref{Thm3} can be proved by the same argument used in Theorem 2.4 in \cite{Li11a}. \hfill $\square$

\section{LSI and lower boundedness of the $\mu$-entropy}

Following \cite{P1}, we introduce the best logarithmic Sobolev constant (called also the $\mu$-entropy)
\begin{eqnarray*}
\mu(\tau):=\inf\limits\left\{W_m(u, \tau): \
u={e^{-f}\over (4\pi \tau)^{m/2}}, \ \int_M ud\mu=1\right\}.
\end{eqnarray*}
In \cite{Li11a}, the author proved that $\mu(\tau)$ is always bounded from below on compact Riemannian manifolds $M$ with $\phi\in C^2(M)$. In this section, we prove a family of logarithmic Sobolev inequalities (LSI) on complete Riemannian manifolds with the uniform volume lower bound condition. The lower boundedness of $\mu(\tau)$ is also proved.

\medskip

\noindent{\bf Proof of Theorem \ref{Thm8a}}. Let $v=\sqrt{u}$.
Then we can rewrite $W_m(u, \tau)$ as follows
\begin{eqnarray*}
W_m(u, \tau)=4\tau\int_M |\nabla v|^2d\mu-\int_M v^2\log
v^2d\mu-\left(m+{m\over 2}\log(4\pi \tau)\right).
\end{eqnarray*}
We need to prove, for all $v\in C_0^\infty(M)$ with $\int_M v^2d\mu=1$, we have
\begin{eqnarray}
4\tau\int_M |\nabla v|^2d\mu-\int_M v^2\log
v^2d\mu-m\left(1+{\over 2}\log(4\pi \tau)\right)\geq \mu(\tau)>-\infty.\label{www}
\end{eqnarray}
Under the conditions of theorem, by Theorem  2.8 in \cite{Li10}, the following $L^2$-Sobolev inequality holds: there exists a constant $C_1(m, K, \mu_0)>0$, depending only on $m, K, \mu_0$,  such
that
\begin{eqnarray*}
\|f\|^2_{2m\over m-2}\leq C_1(m, K, \mu_0)(\|\nabla f||_2^2+\|f\|_2^2), \ \ \ \
\forall f\in C_0^\infty(M).
\end{eqnarray*}
By \cite{Da}, the above $L^2$-Sobolev inequality implies that, for any
$\varepsilon>0$, there exists a constant $\beta(\varepsilon)>0$ such
that the following logarithmic Sobolev inequality holds
\begin{eqnarray}
\int_M f^2\log f^2 d\mu\leq \varepsilon \int_M |\nabla
f|^2+\beta(\varepsilon)\|f\|_2^2+\|f\|_2^2\log\|f\|_2^2,\ \ \ \
\forall f\in C_0^\infty(M),\label{logsob}
\end{eqnarray}
where for some constant $C_2(m, K, \mu_0)>0$ depending only on $m, K, \mu_0$, we have
$$\beta(\varepsilon)\leq C_2(m, K, \mu_0)-{m\over 2}\log\varepsilon.$$
Applying $(\ref{logsob})$ to $f=v\in C_0^\infty(M)$ with $\int_M v^2d\mu=1$, and taking $\varepsilon=4\tau$, we have

\begin{eqnarray*}
\int_M v^2\log v^2 d\mu\leq 4\tau \int_M |\nabla
v|^2-m\left(1+{1\over 2}\log(4\pi\tau)\right)-\mu(\tau),
\end{eqnarray*}
where
$$
-\mu(\tau):=\beta(4\tau)+m\left(1+{1\over 2}\log(4\pi\tau)\right),$$
This proves $(\ref{www})$ with
\begin{eqnarray*}
\mu(\tau)\geq -\left(C_2(m, K, \mu_0)+m+{m\over 2}\log (4\pi)\right), \ \ \
\forall \tau>0.
\end{eqnarray*}
This finishes the proof of Theorem \ref{Thm8a}. \hfill $\square$

\medskip
\noindent{\bf Acknowledgement}. \ \ The author would like to thank F. Baudoin for helpful discussion on the alternative proof of Theorem \ref{Thm0} in September 2011. He would like also to thank D. Bakry, E. Carlen, E. Hsu and A. Thalmaier for their interests on this work.

\medskip

\begin{flushleft}
\medskip\noindent

Xiang-Dong Li, {\sc  Academy of Mathematics and System Science, Chinese
Academy of Sciences, 55, Zhongguancun East Road, Beijing, 100190, P. R. China,}\\
E-mail: xdli@amt.ac.cn
\end{flushleft}

\end{document}